\documentclass[11pt]{amsart}      
\usepackage{amssymb,amsmath}
\usepackage{epsf,latexsym,graphicx}
\usepackage{psfrag}
\title{Minimal Generating sets for the first\\ Syzygies of a Monomial
Ideal}                       
\author{John A. Eagon}       
\date{November 14, 2002}    

\begin{document}             

\maketitle                   

The object of this note is to produce two minimal generating sets of the
first syzygies of a monomial ideal, given the minimal generating set of
the ideal. In the first section we set up notation and give a preliminary result. In the
second section we desribe one small generating set that is not quite minimal. In the third
section we describe a different but analogous small generating set. In the final section we
describe how to reduce both of these to minimal generating sets.

 \section{Notation}

Let $R=k[x_1,x_2,\ldots,x_n]$ be the polynomial ring in $n$ variables over
a field $k$. Let $I$ be a monomial ideal in $R$. Let $G(I)$ be the minimal
generating set of $I$.  For any monomial $m$ in the \emph{lcm-lattice} of
$G(I)$, (i.e.\ the least common multiple lattice of $G(I)$ ordered by
divisibility), define

\begin{center}
         $\Gamma_m = \{\tau: \tau$ is a square free monomial, $\tau|m$, and
	 $\frac{m}{\tau}$ is in $I$\}
\end{center}
where $$\frac{m}{\tau}$$ means the monomial obtained by dividing $m$ by $\tau$.

Also define
\begin{center} $L_{<m} = \{S: S\subset{G(I)}, n_S|m,$ but
$n_S\not=m\}$
\end{center}
where $n_S$ is the least common multiple of the monomials in $S$. Clearly, $n_S$
is in the \emph{lcm-lattice} of $G(I)$. The vertices of $L_{<m}$ may be
identified with the elements of the minimal generating set, $G(I)$.

$\Gamma_m$ is an abstract finite simplicial complex. (Think of the square
free monomials as sets, where divisibility is replaced by containment.) The
facets (maximal faces) of $\Gamma_m$ are among the faces of the form
$$ \sqrt{\frac{m}{\gamma}} $$
where $\gamma$ is in $G(I)$ and $\gamma|m$. Here the root sign
$$ \sqrt{n} $$
stands for the largest square free monomial dividing the monomial $n$.
This is sometimes called the \emph{support} of $n$.

$L_{<m}$ is also an abstract finite simplicial complex. (Use the usual
containment relation.)

For each $\Gamma_m$ which is disconnected, let
$C_{1,m},C_{2.m},\ldots,C_{k_m,m}$ be the connected components (as finite
simplicial complexes).

Select in $C_{i,m}$, for each $i$, one facet of the
form $\sqrt{\frac{m}{\gamma}}$. From all the $\gamma's$ such that
$\sqrt{\frac{m}{\gamma}}$ equals this face select one and label it $\gamma_{i,m}$.
Clearly, we can make this choice because all facets are of this form.

(Later we will make similar, but somewhat more natural choices for $L_{<m}$.)

For each disconnected $\Gamma_m$ let
$$Y_m = \left\{\frac{m}{\gamma_{i,m}}\otimes\gamma_{i,m} -
\frac{m}{\gamma_{j,m}}\otimes\gamma_{j.m}\right\} \subset R\otimes_kR$$
for all pairs $i,j$ such that $1\le i < j\le k_m$. (If $\Gamma_m$ is connected let $Y_m =
\phi$, the empty set.)

\newtheorem{proposition}{Proposition}

\begin{proposition} For each $m$ such that $\Gamma_m$ is disconnected and  for
each pair $\gamma_{i,m},\gamma_{j,m}$ where $i  \not=j$,
$$m=\mathrm{lcm}(\gamma_{i,m},\gamma_{j,m})$$ where the right hand side means
the least common multiple of the two  $\gamma's$. \end{proposition}
\begin{proof}

If not, there exists a variable, $x$, such that
$x|\frac{m}{\gamma_{i,m}}$ and $x|\frac{m}{\gamma_{j,m}}$. But then the faces
$\sqrt{\frac{m}{\gamma_{i,m}}}$ and $\sqrt{\frac{m}{\gamma_{j,m}}}$ overlap and
cannot be in distinct components.

\end{proof}

Now, for all $n$ and $m$ in the \emph{lcm-lattice}, let
 $$Y_{\le{m}} = \bigcup_{n|m} Y_n$$
 $$Y = \bigcup_n\ Y_n$$

For all $\gamma\not=\gamma\prime$ in $G(I)$ let
$n_{\gamma,\gamma\prime}$ be the $lcm(\gamma, \gamma\prime)$.

Consider the set $$\left\{\frac{n_{\gamma, \gamma\prime}}{\gamma}\otimes\gamma -
\frac{n_{\gamma, \gamma\prime}}{\gamma\prime}\otimes\gamma\prime\right\}$$

From Diana Taylor's resolution we conclude that this set generates
the first syzygies of the ideal I.

For a given $m$ in the \emph{lcm-lattice} consider all pairs
$\gamma,\gamma\prime$ such that $n_{\gamma,\gamma\prime} = m$.

Let $$T_m = \left\{\frac{m}{\gamma}\otimes\gamma -
\frac{m}{\gamma\prime}\otimes\gamma\prime\right\}$$ quantified over all such
pairs.

Let $$T = \bigcup_m\ T_m$$

Clearly T is the above set that generates the first syzygies.

If every element of a subset $Q$ of an $R$ module is a linear
combination over $R$ of elements of another subset $P$ we shall say that
``$P$ spans $Q$''. We shall also say that elements of $Q$ are ``in the
span of $P$''.
\section{First Generating Set}
 \newtheorem{theorem}{Theorem}

\begin{theorem} $Y$ spans $T$.

\end{theorem}
\begin{proof}
We will prove by induction on the \emph{lcm-lattice} that
\begin{center}
      $Y_{\le{m}}$ spans $T_m$
\end{center}

Clearly this will suffice to prove the theorem.

Let $\gamma$ and $\gamma\prime$ be any pair of elements in $G(I)$. Let $m
= n_{\gamma,\gamma\prime}$ and consider the element
$$\frac{m}{\gamma}\otimes\gamma - \frac{m}{\gamma\prime}\otimes\gamma\prime$$

It suffices to prove that this is in the span of $Y_{\le{m}}$.

\textbf{Case 1:}

$\sqrt{\frac{m}{\gamma}}$ and $\sqrt{\frac{m}{\gamma\prime}}$ are faces of the
same connected component of $\Gamma_m$.

Then there exists a sequence of elements of $G(I)$, \begin{center} $\gamma =
\gamma_1, \gamma_2, \ldots, \gamma_l = \gamma\prime$ \end{center} such that
$\sqrt{\frac{m}{\gamma_i}}$ and $\sqrt{\frac{m}{\gamma_{i+1}}}$ overlap for $i
= 1, 2, \ldots, l-1$.

So, for each pair $\gamma_i, \gamma_{i+1}$ there exists a variable $x_i$, such
that

   $$x_i |\frac{m}{\gamma_i}\ \text{and}\ x_i |\frac{m}{\gamma_{i+1}}$$

for $i = 1, 2, \ldots, l-1$.
Hence $n_{\gamma_i,\gamma_{i+1}} \not= m$

Let $n_{\gamma_i,\gamma_{i+1}} = n_i,$ for $i = 1, 2,\ldots,l-1$.
Then $n_i|m,$ and $n_i \not=m$.
By induction we assume that $Y_{\le{n_i}}$ spans $T_{n_i}$.

Now, $$\frac{n_i}{\gamma_i}\otimes\gamma_i -
\frac{n_i}{\gamma_{i+1}}\otimes\gamma_{i+1}$$ is in $T_{n_i}$ by the definition
of $T_{n_i}$. We have that  $$Y_{\le{n_i}}\subset{Y_{\le{m}}}$$ since $n_i|m$.

So, $$\frac{n_i}{\gamma_i}\otimes\gamma_i -
\frac{n_i}{\gamma_{i+1}}\otimes\gamma_{i+1}$$  is in the span of $Y_{\le{m}}$
for $i = 1, 2, \ldots, l-1$.

Thus $$\frac{m}{n_i}\left(\frac{n_i}{\gamma_i}\otimes\gamma_i -
\frac{n_i}{\gamma_{i+1}}\otimes\gamma_{i+1}\right)$$
which equals
$$\frac{m}{\gamma_i}\otimes\gamma_i -
\frac{m}{\gamma_{i+1}}
\otimes\gamma_{i+1}$$
is in the span of $Y_{\le{m}}$.

But $$\sum_{i=1}^{l-1}\left(\frac{m}{\gamma_i}\otimes\gamma_i -
\frac{m}{\gamma_{i+1}}\otimes\gamma_{i+1}\right)$$
telescopes to $$\frac{m}{\gamma}\otimes\gamma -
\frac{m}{\gamma\prime}\otimes\gamma\prime$$.
 So this is also in the span of
$Y_{\le{m}}$. Thus Case 1 is proved.

\textbf{Case 2:}

$\sqrt{\frac{m}{\gamma}}$ and $\sqrt{\frac{m}{\gamma\prime}}$ are faces of
different connected components of $\Gamma_m$. Say $$\sqrt{\frac{m}{\gamma}}
\ \text{is in}\ C_{j,m}$$ and $$\sqrt{\frac{m}{\gamma\prime}}\ \text{is in}\
C_{k,m}$$

Thus $\sqrt{\frac{m}{\gamma}}$ and $\sqrt{\frac{m}{\gamma_{j,m}}}$ are in the same
connected component $C_{j,m}$. Hence by Case 1, $$\frac{m}{\gamma}\otimes\gamma -
\frac{m}{\gamma_{j,m}}\otimes\gamma_{j,m}$$ is in the span of $Y_{\le{m}}$.

Similarly, $$\frac{m}{\gamma_{k,m}}\otimes\gamma_{k,m} -
\frac{m}{\gamma\prime}\otimes\gamma\prime$$ is in the span of $Y_{\le{m}}$.

But $$\frac{m}{\gamma_{j,m}}\otimes\gamma_{j,m} -
\frac{m}{\gamma_{k,m}}\otimes\gamma_{k,m}$$ is in $Y_m$ itself.

These three add up to $$\frac{m}{\gamma}\otimes\gamma -
\frac{m}{\gamma\prime}\otimes\gamma\prime$$

\end{proof}
\section {Second Generating Set}

We now give the analogous results for the simplicial complexes $L_{<m}$.

For each $L_{<m}$ which is disconnected let the $C_{i,m}$'s be the connected
components as in the previous case. Let $\gamma_{i,m}$ be any vertex in
$C_{i,m}$. Define the set $Y_{m,L}$ by
$$Y_{m,L} = \left\{\frac{m}{\gamma_{i,m}}\otimes\gamma_{i,m} -
\frac{m}{\gamma_{j,m}}\otimes\gamma_{j.m}\right\} \subset R\otimes_kR$$
using these new $\gamma_{i,m}$'s

\begin{proposition} For each $m$ such that $L_{<m}$ is disconnected and  for
each pair $\gamma_{i,m},\gamma_{j,m}$ where $i\not=j$,
$$m=\mathrm{lcm}(\gamma_{i,m},\gamma_{j,m})$$.
\end{proposition}
\begin{proof}
This is even easier to prove than Proposition 1, since if

     $$m\not=\mathrm{lcm}(\gamma_{i,m},\gamma_{j,m})$$

then there is actually an edge of $L_{<m}$ between the two $\gamma$'s.
\end{proof}

Define
$Y_{\le{m},L}$ and $Y_L$ analogously to $Y_{\le{m}}$ and $Y$

$T_m$, and $T$ are defined as before.

\begin{theorem} $Y_L$ spans $T$. \footnote{The author needed the
first
generating set. He slightly modified the proof given here, which is by Victor Reiner, to
get the proof in the previous section.}
\end{theorem}
\begin{proof}
Exactly as before we prove by induction on the \emph{lcm-lattice} that
\begin{center}
     $Y_{\le{m},L}$ spans $T_m$
\end{center}

Let $\gamma$ and $\gamma\prime$ be any element of $G(I)$.  Let
$m = n_{\gamma,\gamma\prime}$ and consider the element

$$\frac{m}{\gamma}\otimes\gamma - \frac{m}{\gamma\prime}\otimes\gamma\prime$$

It suffices to prove that this is in the span of $Y_{\le{m},L}$.

\textbf{Case 1:}

$\gamma$ and $\gamma\prime$ are vertices of the
same connected component of $L_{<m}$.

Then there exists a sequence of elements of $G(I)$,
\begin{center}

$\gamma = \gamma_1, \gamma_2, \ldots, \gamma_l = \gamma\prime$
\end{center}
such that the edge $\{\gamma_i,\gamma_{i+1}\}$ is in $L_{\le{m}}$ for $i = 1,
2, \ldots, l-1$.

Let $n_i = \mathrm{lcm}(\gamma_i, \gamma_{i+1})$. for $i = 1, 2, \ldots, l-1$.

Then $n_i|m$, but $n_i\not=m$. Thus

$$\frac{n_i}{\gamma}\otimes\gamma_i -
\frac{n_i}{\gamma_{i+1}}\otimes\gamma_{i+1}$$ is in $T_{n_i}$.

By induction we assume that

$$\frac{n_i}{\gamma}\otimes\gamma_i -
\frac{n_i}{\gamma_{I+1}}\otimes\gamma_{i+1}$$

is in the span of $Y_{\le{n_i},L}$ and hence in the span of $Y_{\le{m},L}$.

The argument now goes exactly as before to complete Case 1.

\textbf{Case 2:}

$\gamma$ and $\gamma\prime$ are vertices of different connected components of
$L_{\le{m}}$.

The argument again goes as before, and I leave it as an exercise to complete
the proof of the theorem.
\end{proof}

\section{Minimal Generating Sets for the\\ First Syzygies}

By Theorem 2.1, page 523 of \cite{GPW}, and Proposition 1.1 of \cite{BH} we
have that the minimal number of 1st syzygies of $I$ (which is the same as the
minimal number of 2nd syzygies of $R/I$) is given by the formula

$$b_2(R/I) = \sum_m\ \dim\tilde{H}_{0}(\Gamma_m) = \sum_m\
\dim\tilde{H}_{0}(L_{<m})$$  where $m$ is quantified over all monomials in
the \emph{lcm-lattice}

For each $m$ in the \emph{lcm-lattice}, $\dim\tilde{H}_0(\Gamma_m) =
\dim\tilde{H}_0(L_{<m})$ is one less than $k_m$, the number of connected
components of $\Gamma_m$, respectively, $L_{<m}$.

For each $m$ such that $\Gamma_m$, respectively, $L_{<m}$ is disconnected we
consider the set  $$Z_m,\ \text{respectively,}\ Z_{m,L} = \left\{\frac{m}{\gamma_{i,m}}\otimes\gamma_{i,m} -
\frac{m}{\gamma_{1,m}}\otimes\gamma_{1,m}\right\}$$ for $i = 2, 3, \ldots,
k_m$.

The cardinality of this set is also one less than $k_m$. On the other
hand, by taking differences, we see that this set spans $Y_m$,
respectively,$Y_{m,L}$.
Quantifying over all $m$ in the \emph{lcm-lattice}, we see that the set
  $$Z,\ \text{respectively,}\ Z_L = \bigcup_m\ Z_m,\ \text{respectively,}\
  Z_{m,L}$$
spans $Y$, respectively $Y_L$ and hence $T$ by Theorem 1, respectively Theorem
2.

By the formula cited above from \cite{GPW} we see that the cardinality  of $Z$,
respectively, $Z_L$ is the number of minimal generators of the first syzygies.

Since $Z$, respectively $Z_L$ has the correct cardinality and spans $T$, each
must be a  minimal generating set for the first syzygies.


\begin{thebibliography}{GPW}
\bibitem[BH]{BH}W. Bruns and J. Herzog, \it Semigroup rings and simplicial
complexes, \rm J. Pure. Appl. Algebra \bf 122 \rm (1997), 185-208.
\bibitem[GPW]{GPW}V. Gasharov, I. Peeva and V. Welker, \it The lcm-lattice
in monomial resolutions, \rm Math. Res. Lett. \bf 6, \rm (1999), 521-532.
\end{thebibliography}
\end{document}